\documentclass[a4paper,10pt,twoside]{article}

\usepackage{amssymb}
\usepackage{amsthm}
\usepackage{amsmath}
\usepackage{enumerate}
\usepackage{hyperref}

\newtheorem{thmA}{Theorem}
\newtheorem{thm}{Theorem}
\newtheorem{lemA}{Lemma}
\newtheorem{lem}{Lemma}

\newtheorem{cor}{Corollary}

\newcommand*{\half}{\frac12}
\newcommand*{\quart}{\frac14}
\newcommand*{\lb}{\left\{}
\newcommand*{\rb}{\right\}}
\newcommand*{\PP}{\mathbb{P}}
\newcommand*{\di}{\, \mathrm{d} }
\newcommand*{\qv}[2][M]{[ #1 ]_{#2}}
\newcommand*{\loc}{\ensuremath{\mathcal{L}}}
\newcommand*{\ind}{\ensuremath{\mathbf{1}_}}
\newcommand*{\sign}{\ensuremath{\mathrm{sgn}}}

\begin{document}

\title{Stochastic integration based on simple, symmetric
random walks}

\author{Tam\'as Szabados\footnote{Corresponding author, address:
Department of Mathematics, Budapest University of Technology and
Economics, M\H{u}egyetem rkp. 3, H \'ep. V em. Budapest, 1521,
Hungary, e-mail: szabados@math.bme.hu, telephone: (+36 1)
463-1111/ext. 5907, fax: (+36 1) 463-1677} \footnote{Research
supported by a Hungarian National Research Foundation (OTKA) grant
No. T42496.} and Bal\'azs Sz\'ekely \footnote{Research
supported by the HSN laboratory of BUTE.} \\
Budapest University of Technology and Economics}

\date{}

\maketitle

\bigskip

%%%%%%%%%%%%%%%%%%%%%%%%%%%%%%%%%%%%%%%%%%%%%%%%%%%%%%%%%%%%%%%%%%%

\begin{abstract}

A new approach to stochastic integration is described, which is
based on an a.s. pathwise approximation of the integrator by
simple, symmetric random walks. Hopefully, this method is
didactically more advantageous, more transparent, and technically
less demanding than other existing ones. In a large part of the
theory one has a.s. uniform convergence on compacts. In
particular, it gives a.s. convergence for the stochastic integral
of a finite variation function of the integrator, which is not
c\`adl\`ag in general.

\end{abstract}

%%%%%%%%%%%%%%%%%%%%%%%%%%%%%%%%%%%%%%%%%%%%%%%%%%%%%%%%%%%%%%%%%%%

\renewcommand{\thefootnote}{\alph{footnote}}
\footnotetext{ 2000 \emph{MSC.} Primary 60H05. Secondary 60F15,
60G50.} 

%%%%%%%%%%%%%%%%%%%%%%%%%%%%%%%%%%%%%%%%%%%%%%%%%%%%%%%%%%%%%%%%%%%

\section{Introduction} \label{sec:Intro}

The main purpose of the present paper is to describe a new
approach to stochastic integration, which is based on an a.s.
pathwise approximation of the integrator by simple, symmetric
random walks. It is hoped that this method is pedagogically more
advantageous, more transparent, and technically less demanding
than other existing ones. This way hopefully a larger,
mathematically less mature audience may get acquainted with a
theory of stochastic integration. Since stochastic calculus plays
an ever-increasing role in several applications (mathematical
finance, statistical mechanics, engineering, \dots) nowadays, this
aim seems to be justified.

A basic feature of the present theory that, except for the
integral of general predictable integrands, one always has almost
sure uniform convergence on compact intervals. This is true for
the approximation of the integrators, quadratic variations, local
times, and stochastic integrals when the integrand is a $C^1$ or a
finite variation function of the integrator. We believe that for a
beginner it is easier to understand almost sure pathwise
convergence than $L^2$ or in probability convergence that
typically appear in stochastic calculus textbooks. We mention that
our method significantly differs from the a.s. converging
approximation given by Karandikar \cite{Kar1995}. Important tools
in our approach are discrete versions of It\^o's and It\^o--Tanaka
formulas. The continuous versions easily follow from these by a.s.
pathwise convergence. (See earlier versions of it in
\cite{Szab1990} and  \cite{Szab1996}.) To our best knowledge, our
approach is new in giving a.s. convergence for the stochastic
integral of a finite variation function of the integrator, which
is not c\`adl\`ag in general.

In the case of a more general, e.g. predictable integrand, our
method is closer to the traditional ones, e.g. we too have $L^2$
convergence then. The only, hopefully interesting, difference is
that in the approximating sums random partitions are used so that
the values of the integrand are multiplied by $\pm 1$ times a
constant scaling factor. The signs are from an independent,
symmetric coin tossing sequence.

The most general integrators considered in this paper are
continuous local martingales $M$. It is easy to extend this to
semimartingales $X$ that can be written $X = M + A$, where $M$ is
a continuous local martingale and $A$ is a finite variation
process. The reason is simple: integration with respect to such an
$A$ can also be defined pathwise.

%%%%%%%%%%%%%%%%%%%%%%%%%%%%%%%%%%%%%%%%%%%%%%%%%%%%%%%%%%%%%%%%%%%

\section{Preliminaries} \label{sec:Pre}

A basic tool of the present paper is an elementary construction of
Brownian motion (BM). The specific construction we are going to
use in the sequel, taken from \cite{Szab1996}, is based on a
nested sequence of simple, symmetric random walks (RW's) that
uniformly converges to the Wiener process (=BM) on bounded
intervals with probability $1$. This will be called \emph{``twist
and shrink''} construction in the sequel. Our method is a
modification of the one given by Frank Knight in 1962
\cite{Kni62}.

We summarize the major steps of the ``twist and shrink''
construction here. We start with an infinite matrix of independent
and identically distributed random variables $X_m(k)$, $\PP \lb
X_m(k)= \pm 1 \rb = \frac12$ ($m\ge 0$, $k\ge 1$), defined on the
same complete probability space $(\Omega,\mathcal{F},\PP)$. (All
stochastic processes in the sequel will be defined on this
probability space.) Each row of this matrix is a basis of an
approximation of the Wiener process with a dyadic step size
$\Delta t=2^{-2m}$ in time and a corresponding step size $\Delta
x=2^{-m}$ in space. Thus we start with a sequence of independent
simple, symmetric RW's $S_m(0) = 0$, $S_m(n) = \sum_{k=1}^{n}
X_m(k)$ $(n \ge 1)$.

The second step of the construction is \emph{twisting}. From the
independent RW's we want to create dependent ones so that after
shrinking temporal and spatial step sizes, each consecutive RW
becomes a refinement of the previous one.  Since the spatial unit
will be halved at each consecutive row, we define stopping times
by $T_m(0)=0$, and for $k\ge 0$,
\[
T_m(k+1)=\min \{n: n>T_m(k), |S_m(n)-S_m(T_m(k))|=2\} \qquad (m\ge
1)
\]
These are the random time instants when a RW visits even integers,
different from the previous one. After shrinking the spatial unit
by half, a suitable modification of this RW will visit the same
integers in the same order as the previous RW. We operate here on
each point $\omega\in\Omega$ of the sample space separately, i.e.
we fix a sample path of each RW. We define twisted RW's
$\tilde{S}_m$ recursively for $k=1,2,\dots$ using
$\tilde{S}_{m-1}$, starting with $\tilde{S}_0(n)=S_0(n)$ $(n\ge
0)$ and $\tilde{S}_m(0) = 0$ for any $m \ge 0$. With each fixed
$m$ we proceed for $k=0,1,2,\dots$ successively, and for every $n$
in the corresponding bridge, $T_m(k)<n\le T_m(k+1)$. Any bridge is
flipped if its sign differs from the desired:
\[
\tilde{X}_m(n)=\left\{
\begin{array}{rl}
 X_m(n)& \mbox{ if } S_m(T_m(k+1)) - S_m(T_m(k))
= 2\tilde X_{m-1}(k+1), \\
- X_m(n)& \mbox{ otherwise,}
\end{array}
\right.
\]
and then $\tilde{S}_m(n)=\tilde{S}_m(n-1)+\tilde{X}_m(n)$. Then
$\tilde{S}_m(n)$ $(n\ge 0)$ is still a simple symmetric RW
\cite[Lemma 1]{Szab1996}. The twisted RW's have the desired
refinement property:
\[
\tilde{S}_{m+1}(T_{m+1}(k)) = 2 \tilde{S}_{m}(k) \qquad (m\ge 0,
k\ge 0).
\]

The third step of the RW construction is \emph{shrinking}. The
sample paths of $\tilde{S}_m(n)$ $(n\ge 0)$ can be extended to
continuous functions by linear interpolation, this way one gets
$\tilde{S}_m(t)$ $(t\ge 0)$ for real $t$. The $mth$ \emph{``twist
and shrink'' RW} is defined by
\[
\tilde{B}_m(t)=2^{-m}\tilde{S}_m(t2^{2m}).
\]
Then the \emph{refinement property} takes the form
\begin{equation}
\tilde{B}_{m+1}\left(T_{m+1}(k)2^{-2(m+1)}\right) = \tilde{B}_m
\left( k2^{-2m}\right) \qquad (m\ge 0,k\ge 0). \label{eq:refin}
\end{equation}
Note that a refinement takes the same dyadic values in the same
order as the previous shrunken walk, but there is a \emph{time
lag} in general:
\begin{equation} T_{m+1}(k)2^{-2(m+1)} - k2^{-2m} \ne 0 .
\label{eq:tlag}
\end{equation}

It is clear that this construction is especially useful for local
times, since a refinement approximates the local time of the
previous walk, with a geometrically distributed random number of
visits with half-length steps (cf. Theorem \ref{th:loctime+}
below).

Now we quote some important facts from \cite{Szab1996} and
\cite{SzaSze2005} about the above RW construction that will be
used in the sequel.

\begin{thmA} \label{th:Wiener}
On bounded intervals the sequence $(\tilde B_m)$ almost surely
uniformly converges as $m \to \infty$ and the limit process is
Brownian motion $W$. For any $C>1$, and for any $K>0$ and $m \ge
1$ such that $K 2^{2m} \ge N(C)$, we have
\begin{eqnarray*}
\PP \lb \sup_{0 \le t \le K} |W(t) - \tilde{B}_m(t)| \ge 27 \: C
K_*^{\quart} (\log_*K)^{\frac34} m^{\frac34} 2^{-\frac{m}{2}} \rb
\\
\le \frac{6}{1-4^{1-C}} (K2^{2m})^{1-C} ,
\end{eqnarray*}
where $K_* := K \vee 1$ and $\log_* K := (\log K) \vee 1$.

\end{thmA}
($N(C)$ here and in the sequel denotes a large enough integer
depending on $C$, whose value can be different at each occasion.)

Conversely, with a given Wiener process $W$, one can define the
stopping times which yield the \emph{Skorohod embedded RW's}
$B_m(k2^{-2m})$ into $W$. For every $m\ge 0$ let $s_m(0)=0$ and
\begin{equation} \label{eq:Skor1}
s_m(k+1)=\inf{}\{s: s > s_m(k), |W(s)-W(s_m(k))|=2^{-m}\} \qquad
(k \ge 0).
\end{equation}
With these stopping times the embedded dyadic walks by definition
are
\begin{equation} \label{eq:Skor2}
B_m(k2^{-2m}) = W(s_m(k)) \qquad (m\ge 0, k\ge 0).
\end{equation}
This definition of $B_m$ can be extended to any real $t \ge 0$ by
pathwise linear interpolation.

If a Wiener process is built by the ``twist and shrink''
construction described above using a sequence $(\tilde{B}_m)$ of
nested RW's and then one constructs the Skorohod embedded RW's
$(B_m)$, it is natural to ask about their relationship. The next
lemma explains that they are asymptotically identical. In general,
roughly saying, $(\tilde{B}_m)$ is more useful when someone wants
to generate stochastic processes from scratch, while $(B_m)$ is
more advantageous when someone needs discrete approximations of
given processes.

\begin{lemA} \label{le:equid}
For any $C > 1$, and for any $K>0$ and $m \ge 1$ such that $K
2^{2m} \ge N(C)$ take the following subset of $\Omega $:
\[
 A^*_{K,m} = \lb \sup_{n > m} \: \sup_{k}
|2^{-2n} T_{m,n}(k) - k2^{-2m}| < (42 \: C K \log_*K)^{\half}
m^{\half} 2^{-m} \rb ,
\]
where $T_{m,n}(k) = T_n \circ T_{n-1} \circ \cdots \circ T_m(k)$
for $n>m\ge 0$ and $k \in [0, K2^{2m}]$. Then
\[
\PP \lb (A^*_{K,m})^c \rb \le \frac{2}{1 - 4^{1-C}}(K2^{2m})^{1-C}
.
\]
Moreover, $\lim_{n\to \infty} 2^{-2n} T_{m,n}(k) = t_m(k)$ exists
almost surely and on $A^*_{K,m}$ we have
\[
\tilde{B}_m(k2^{-2m}) = W(t_m(k)) \qquad (0 \le k2^{-2m} \le K) ,
\]
cf. (\ref{eq:Skor2}). Further, almost everywhere on $A^*_{K,m}$,
$s_m(k) = t_m(k)$ and
\[
\sup_{0\le k2^{-2m}\le K} |s_m(k)-k2^{-2m}| \le (42 \: C K \log_*
K) ^{\half} m^{\half} 2^{-m} .
\]
\end{lemA}

The next theorem shows that the rate of convergence of $(B_m)$ to
$W$ is essentially the same as the one of $(\tilde{B}_m)$, cf.
Theorem \ref{th:Wiener}.

\begin{thmA} \label{th:Wiener_Skor}
For any $C > 1$, and for any $K>0$ and $m \ge 1$ such that $K
2^{2m} \ge N(C)$ we have
\begin{eqnarray*}
\PP \lb \sup_{0\le t\le K} \left| W(t) - B_m(t) \right| \ge 27 \:
C K_*^{\quart} (\log_*K)^{\frac34} m^{\frac34}
2^{-\frac{m}{2}} \rb \\
\le \frac{8}{1 - 4^{1-C}} (K2^{2m})^{1-C} .
\end{eqnarray*}

\end{thmA}

In the last part of this section we quote a result from
\cite{SzaSze2005} about an elementary definition of Brownian local
time, based on the ``twist and shrink'' RW's $(\tilde{B}_m(t))$.
This is basically a finite-time-horizon version of a theorem of
R\'ev\'esz \cite{Rev81}, in a somewhat sharper form.

First, one can define one-sided, \emph{up and down local times}
$\tilde{\ell}^+_m(k,x)$ and $\tilde{\ell}^-_m(k,x)$ of the
``twisted'' RW's $(\tilde{S}_m(j))$ at a point $x \in \mathbb{Z}$
up to time $k \in \mathbb{Z}_+$ as $\tilde{\ell}^{\pm}_m(0,x) = 0$
and
\begin{equation} \label{eq:loctime_int}
\tilde{\ell}^{\pm}_m(k,x) = \#\{j : 0 \le j < k, \: \tilde{S}_m(j)
= x, \: \tilde{S}_m(j+1) = x \pm 1 \}  \quad  (k \ge 1).
\end{equation}
Then the two-sided local time is $\tilde{\ell}_m(k,x) =
\tilde{\ell}^+_m(k,x) + \tilde{\ell}^-_m(k,x)$.

The one-sided (and two-sided) local times of the $m$th ``twist and
shrink'' walk $\tilde B_m$ at a point $x \in 2^{-m} \mathbb{Z}$ up
to time $t \in 2^{-2m} \mathbb{Z}_+$ are defined as
\begin{equation} \label{eq:loctime}
\tilde{\loc}_m^{\pm}(t,x) = 2^{-m} \tilde{\ell}_m^{\pm}
\left(t2^{2m}, x 2^m \right),
\end{equation}
corresponding to the fact that the spatial step size of $\tilde
B_m$ is $2^{-m}$. Then $\tilde{\loc}_m^{\pm}(t,x)$ can be extended
to any $t \in \mathbb{R}_+$ and $x \in \mathbb{R}$ by linear
interpolation, as a continuous two-parameter process.

\begin{thmA} \label{th:loctime+}
On any strip $[0, K] \times \mathbb{R}$ the sequence
$(\tilde{\loc}^+_m(t,x))$ almost surely uniformly converges as $m
\to \infty$ to the one half of the Brownian local time
$\loc(t,x)$. For any $C>1$, and for any $K>0$ and $m \ge 1$ such
that $K 2^{2m} \ge N(C)$, we have
\begin{eqnarray*}
\PP \lb \sup_{(t,x) \in [0,K] \times \mathbb{R}} \left|\half \:
\loc(t,x) - \tilde{\loc}^+_{m+1}(t,x)\right| \ge 50 \: C
K_*^{\quart}
(\log_*K)^{\frac34} m^{\frac34} 2^{-\frac{m}{2}} \rb \\
\le \frac{30}{1-4^{1-C}} (K2^{2m})^{1-C} .
\end{eqnarray*}
Similar statements hold for $(\tilde{\loc}^-_m(t,x))$ as well.
\end{thmA}

Very similar statements hold when the ``twist and shrink'' walks
$\tilde{B}_m(t)$ are replaced by Skorohod embedded walks $B_m(t)$
in the definition (\ref{eq:loctime_int}) and (\ref{eq:loctime}) of
local time. This local time will be denoted by $\loc^{\pm}_m(t,x)$
in the sequel.

\begin{cor} \label{co:loctime+}
On any strip $[0, K] \times \mathbb{R}$ the sequence
$(\loc^+_m(t,x))$ almost surely uniformly converges as $m \to
\infty$ to the one half of the Brownian local time $\loc(t,x)$.
For any $C>1$, and for any $K>0$ and $m \ge 1$ such that $K 2^{2m}
\ge N(C)$, we have
\begin{eqnarray*}
\PP \lb \sup_{(t,x) \in [0,K] \times \mathbb{R}} \left|\half \:
\loc(t,x) - \loc^+_{m+1}(t,x)\right| \ge 50 \: C K_*^{\quart}
(\log_*K)^{\frac34} m^{\frac34} 2^{-\frac{m}{2}} \rb \\
\le \frac{32}{1-4^{1-C}} (K2^{2m})^{1-C} .
\end{eqnarray*}
Similar statements hold for $(\loc^-_m(t,x))$ as well.
\end{cor}

\begin{proof}
By Lemma \ref{le:equid}, almost everywhere on an event $A^*_{K,m}$
whose complement has probability going to zero exponentially fast
with $m$, one has that
\[
\tilde{B}_m(k 2^{-2m}) = W(t_m(k)) = W(s_m(k)) = B_m(k 2^{-2m})
\qquad (0 \le k2^{-2m} \le K).
\]
Thus on $A^*_{K,m}$, for almost every $\omega$ and for all $(t,x)
\in [0,K] \times \mathbb{R}$, $\tilde{\loc}^+_m(t,x) =
\loc^+_m(t,x)$.

Hence, by the triangle inequality
\[
\sup_{(t,x)} |\loc^+(t,x) - \loc^+_m(t,x)| \le \sup_{(t,x)}
|\tilde{\loc}^+_m(t,x) - \loc^+_m(t,x)| + \sup_{(t,x)}
|\loc^+(t,x) - \tilde{\loc}^+_m(t,x)|,
\]
and by Theorem \ref{th:loctime+}, we obtain the above statement:
\begin{eqnarray*}
\lefteqn{\PP \lb \sup_{(t,x) \in [0,K] \times \mathbb{R}}
\left|\half \: \loc(t,x) - \loc^+_{m+1}(t,x)\right| \ge 50 \: C
K_*^{\quart} (\log_*K)^{\frac34} m^{\frac34} 2^{-\frac{m}{2}}
\rb} \\
&\le& \PP \lb (A^*_{K,m})^c \rb \\
&+& \PP \lb \sup_{(t,x) \in [0,K] \times \mathbb{R}} \left|\half
\: \loc(t,x) - \tilde{\loc}^+_{m+1}(t,x)\right| \ge 50 \: C
K_*^{\quart}
(\log_*K)^{\frac34} m^{\frac34} 2^{-\frac{m}{2}} \rb \\
&\le& \frac{32}{1-4^{1-C}} (K2^{2m})^{1-C} .
\end{eqnarray*}

\end{proof}

%%%%%%%%%%%%%%%%%%%%%%%%%%%%%%%%%%%%%%%%%%%%%%%%%%%%%%%%%%%%%%%%%%%

\section{A discrete It\^o--Tanaka formula and its consequences}
\label{sec:Ito}

It is interesting that one can give discrete versions of It\^o's
formula and of It\^o--Tanaka formula, which are purely algebraic
identities, not assigning any probabilities to the terms. Despite
this, the usual It\^o's and It\^o--Tanaka formulae follow fairly
easily in a proper probability setting.

Discrete It\^o formulas are not new. Apparently, the first such
formula was given by Kudzma in 1982 \cite{Kud1982}. A similar
approach was recently used by Fujita \cite{Fuj2003} and Akahori
\cite{Aka2005}. Discrete Tanaka formulae were given by Kudzma
\cite{Kud1982} and Cs\"org\H{o}--R\'ev\'esz \cite{CsoRev1985}. The
elementary algebraic approach used in the present paper is
different from these; it was introduced by the first author in
1989 \cite{Szab1990}.

First consider an arbitrary numerical sequence $X_r = \pm 1$ ($r
\ge 1$) and partial sums $S_0=a \in \mathbb{R}$, $S_n =
a+(X_1+\dots +X_n) \Delta x$ ($n \ge 1$), where $\Delta x > 0$ is
the step-size. Let $f:\mathbb{R} \rightarrow \mathbb{R}$, and $b =
a + k \Delta x$, $k \in \mathbb{Z}$. We define \emph{trapezoidal
sums} of $f$ with step $\Delta x$ on the interval $[a, b]$ by
\begin{equation*}
T_{x=a}^b f(x) \Delta x := \sign (b-a) \Delta x \left\{\frac{1}{2}
f(a) + \sum_{j=1}^{|k|-1} f(a+\sign (b-a) j\Delta x ) +
\frac{1}{2} f(b) \right\}
\end{equation*}
where $\sign(x) = 1, 0 , -1$ according as $x > 0$, $x = 0$, or $x
< 0$, respectively.

Similarly as above, one-sided, \emph{up and down local times}
$\loc^+(n,x)$ and $\loc^-(n,x)$ ($n \ge 0$) of the series $(S_n)$
can be defined with step $\Delta x$: $\loc^{\pm}(0,x) := 0$ and
\[
\loc^{\pm}(n,x) := \Delta x \: \#\{j : 0 \le j < n, \: S_{j} = x,
\: S_{j+1} = x \pm \Delta x \} ,
\]
where $n \ge 1$ and $x = a + k \Delta x$, $k \in \mathbb{Z}$. The
two-sided local time is $\loc(n,x) := \loc^+(n,x) + \loc^-(n,x)$.

\begin{lem}\label{le:discrete_Ito}
Take any function $f:\mathbb{R}\rightarrow\mathbb{R}$ , $a \in
\mathbb{R}$, step $\Delta x > 0$, and numerical sequence $X_r =
\pm 1$ ($r \ge 1$). Let $S_0=a$, $S_n=a+(X_1+\dots +X_n)\Delta x$
($n \ge 1$). Then the following equalities hold:
\begin{equation} \label{eq:disc_Strat}
T_{x=S_0}^{S_n} f(x) \Delta x = \sum_{r=1}^{n} \frac{f(S_r) +
f(S_{r-1})}{2} X_r \Delta x
\end{equation}
(discrete Stratonovich formula). Alternatively,
\begin{equation} \label{eq:disc_Ito}
T_{x=S_0}^{S_n} f(x) \Delta x = \sum_{r=1}^{n} f(S_{r-1}) X_r
\Delta x + \frac12 \sum_{r=1}^n \frac{f(S_r) - f(S_{r-1})}{X_r
\Delta x} (\Delta x)^2
\end{equation}
(discrete It\^o's formula). Finally,
\begin{eqnarray} \label{eq:disc_Ito_Tan}
\lefteqn{ T_{x=S_0}^{S_n} f(x)\Delta x = \sum_{r=1}^{n}
f(S_{r-1}) X_r\Delta x}  \\
&+& \frac12 \sum_{x \in a+\mathbb{Z} \Delta x} \left\{\loc^+(n,x)
(f(x+\Delta x)  - f(x)) + \loc^-(n,x) (f(x) - f(x-\Delta x))
\right\}  \nonumber
\end{eqnarray}
(discrete It\^o--Tanaka formula).
\end{lem}

\begin{proof}
Algebraically,
\begin{eqnarray*}
\lefteqn{T_{x=S_0}^{S_r} f(x) \Delta x - T_{x=S_0}^{S_{r-1}} f(x)
\Delta x = X_r \Delta x \frac{f(S_{r-1})+f(S_r)}{2}} \\
& = & f(S_{r-1})X_r \Delta x + \frac12 \frac{f(S_r)-f(S_{r-1})}
{X_r \Delta x} (\Delta x)^2 \\
& = & f(S_{r-1}) X_r \Delta x \\
&&+ \frac12 \Delta x \sum_{x \in a+\mathbb{Z} \Delta x}
(f(x+\Delta x) - f(x))
\ind{\{S_{r-1}=x, S_r=x+\Delta x\}} \\
&&+ \frac12 \Delta x \sum_{x \in a+\mathbb{Z} \Delta x}
(f(x)-f(x-\Delta x)) \ind{\{S_{r-1}=x, S_r=x-\Delta x\}}.
\end{eqnarray*}
The first equality follows from the fact that if $X_r=1$, one has
to add a term, while if $X_r=-1$, one has to subtract a term in
the trapezoidal sum. Then the second equality follows since $1/X_r
= X_r$, and the third equality is obviously implied by the second.
Summing up for $r=1, \dots, n$, the sum on the left telescopes,
and from the three equalities one obtains the three formulae,
respectively.
\end{proof}

Introducing the notation $h_{\pm \Delta x}(x) := (f(x \pm \Delta
x) - f(x))/(\pm \Delta x)$ and comparing (\ref{eq:disc_Ito}) and
(\ref{eq:disc_Ito_Tan}), we obtain a \emph{discrete occupation
time formula}, cf. (\ref{eq:occup}):
\[
\sum_{r=1}^n h_{X_r \Delta x}(S_{r-1})   (\Delta x)^2 = \sum_{x
\in a+\mathbb{Z} \Delta x} \left\{h_{\Delta x}(x) \loc^+(n,x) +
h_{-\Delta x}(x) \loc^-(n,x)  \right\} \Delta x .
\]

Let us apply Lemma \ref{le:discrete_Ito} with Skorohod embedded
walks $(B_m(t))$. Elementary calculus shows, c.f. \cite[Theorem
6]{Szab1990}, that when $g \in C^2(\mathbb{R})$ and we set $f =
g'$ in (\ref{eq:disc_Strat}) or in (\ref{eq:disc_Ito}), then, as
$m \to \infty$, the terms converge almost surely to the
corresponding terms of the Stratonovich and It\^o's formula,
respectively. On one hand, this gives a definition of the It\^o
integral $\int_{0}^{t} f(W(s)) \di W(s)$ and Stratonovich integral
$\int_{0}^{t} f(W(s)) \: \partial W(s)$, with almost sure
convergence of It\^o sums
\begin{equation} \label{eq:stocsum1}
\left(f(W) \cdot W \right)^m_t := \sum_{r=1}^{\lfloor t 2^{2m}
\rfloor} f(B_m((r-1) 2^{-2m})) \: 2^{-m} X_m(r),
\end{equation}
and Stratonovich sums
\[
\left(f(W) \circ W \right)^m_t := \sum_{r=1}^{\lfloor t 2^{2m}
\rfloor} \frac{f(B_m((r-1) 2^{-2m})) + f(B_m(r 2^{-2m}))}{2} \:
2^{-m} X_m(r),
\]
respectively. Here $(X_m(r))_{r=1}^{\infty}$ is an independent,
$\pm 1$ symmetric coin tossing sequence, $2^{-m} X_m(r) = B_m(r 2
^{-2m}) - B_m((r-1) 2^{-2m})$. This essentially means that in the
sums approximating the It\^o integral we apply a partition with
the Skorohod stopping times $0 = s_m(0) < s_m(1) < \cdots <
s_m(\lfloor t 2^{2m} \rfloor)$, since $B_m(r2^{-2m}) = W(s_m(r))$.
Second, taking almost sure limits as $m \to \infty$, one
immediately obtains the corresponding It\^o's formula
\begin{equation} \label{eq:Ito_form}
g(W(t)) - g(W(0)) = \int_{0}^{t} g'(W(s)) \di W(s) + \frac12
\int_{0}^{t} g''(W(s)) \di s ,
\end{equation}
and Stratonovich formula
\[
g(W(t)) - g(W(0)) = \int_{0}^{t} g'(W(s)) \: \partial W(s) .
\]

In the same way, now we show almost sure convergence of stochastic
sums $(f(W) \cdot W)_{t}^{m}$, when $g$ is the difference of two
convex functions and $f=g'_-$, its left derivative. Then we
immediately obtain the It\^o--Tanaka formula as well, with the
help of (\ref{eq:disc_Ito_Tan}).

\begin{thm} \label{th:non-cadlag_int}
Let $g$ be the difference of two convex functions, $g'_-$ be its
left derivative, and let $\mu$ be the signed measure defined by
the change of $g'_-$ when restricted to compacts (the second
derivative of $g$ in the generalized function sense). Then for
arbitrary $K>0$,
\begin{equation} \label{eq:non-cadlag_int}
\sup_{t\in[0,K]}\left|\left(g'_- (W) \cdot W \right)_t^m -
\int_0^t g'_- (W(s)) \di W(s) \right| \rightarrow 0
\end{equation}
almost surely as $m \to \infty$, and for any $t \ge 0$,
\begin{equation} \label{eq:Ito_Tan_form}
g(W(t)) - g(W(0)) = \int_{0}^{t} g'_- (W(s)) \di W(s) + \frac12
\int_{-\infty}^{\infty} \loc (t,x) \mu(\di x)  .
\end{equation}
\end{thm}

\begin{proof}
The basic idea of the proof is that one may substitute Skorohod
embedded walks $B_m(r2^{-2m}) = W(s_m(r))$, $B_m(0) = W(0) = a \in
\mathbb{R}$ into (\ref{eq:disc_Ito_Tan}) to obtain
\begin{eqnarray} \label{eq:disc_Ito_Tan_Skor}
\left(g'_- (W) \cdot W \right)_t^m &=& T_{x=a}^{B_m(\lfloor t
2^{2m} \rfloor 2^{-2m})} g'_-(x) 2^{-m}
\nonumber \\
&-& \frac12 \sum_{x \in a + 2^{-m} \mathbb{Z}}
\left\{\loc^+_m(t,x) (g'_-(x+2^{-m}) - g'_-(x)) \right. \nonumber
\\
&& \qquad \qquad \left. + \loc^-_m(t,x) (g'_-(x) - g'_-(x-2^{-m}))
\right\} .
\end{eqnarray}
Then it is enough to deal with the almost sure uniform convergence
on compacts of the terms on the right side of
(\ref{eq:disc_Ito_Tan_Skor}), and that will imply the same
convergence of the stochastic sum on the left side and will prove
the It\^o--Tanaka formula (\ref{eq:Ito_Tan_form}) as well.

First we show (\ref{eq:Ito_Tan_form}) for measures $\mu$ supported
in a compact interval $[-M,M]$. Because of the linearity of
(\ref{eq:Ito_Tan_form}), it can be supposed that $g$ is convex.
Then $g'_-$ is non-decreasing, left-continuous, and $g'_-(x)$ is
constant if $x > M$ or $x \le -M$, hence $g'_-$ is bounded on
$\mathbb{R}$. The measure $\mu$, defined by $\mu([a,b)) := g'_-(b)
- g'_-(a)$ if $a < b$, is a regular, finite, positive Borel
measure on $\mathbb{R}$, with total mass $\mu(\mathbb{R}) =
g'_+(M) - g'_-(-M)$.

We are going to prove (\ref{eq:Ito_Tan_form}) pathwise. For this,
let $\Omega_0$, $\PP \lb \Omega_0 \rb = 1$, denote a subset of the
sample space $\Omega $, on which, as $m \to \infty$, $B_m(t)$
uniformly converges to $W(t)$ on $[0,K]$ and $\loc^{\pm}_m(t,x)$
uniformly converges to $\loc^{\pm}(t,x)$ on $[0,K]\times
\mathbb{R}$, cf. Theorem \ref{th:Wiener_Skor} and Corollary
\ref{co:loctime+} above. We fix an $\omega \in \Omega _0$. Then,
obviously, $W(t)$ and $\loc^{\pm}(t,x)$ will have continuous
paths.

Consider the first term on the right side of
(\ref{eq:disc_Ito_Tan_Skor}). We want to show that it uniformly
converges to $g(W(t)) - g(W(0))$ for $t \in [0,K]$. With the
notation $t_m:=\lfloor t 2^{2m} \rfloor 2^{-2m}$,
\begin{eqnarray} \label{eq:trapez}
\lefteqn{\sup_{t\in [0,K]} \left|T_{x=a}^{B_m(t_m)} g_-'(x) 2^{-m}
- \int_{a}^{W(t)} g_-'(x) \di x \right|} \nonumber \\
&\le& \sup_{t\in [0,K]} \left|T_{x=a}^{B_m(t_m)} g_-'(x) 2^{-m} -
\int_a^{B_m(t_m)} g_-'(x) \di x \right| \nonumber \\
&+& \sup_{t\in [0,K]} \left| \int_{B_m(t_m)}^{W(t)} g_-'(x) \di x
\right|.
\end{eqnarray}

The first term on the right side of (\ref{eq:trapez}) goes to zero
as $m \to \infty$, because $g'_-$ is non-decreasing, hence Riemann
integrable on any compact interval $[a,B_m(t_m)]$, the trapezoidal
sum is a Riemann sum on the same interval, so their difference can
be estimated from above by the difference of the upper sum and the
lower sum of $g'_-$ on $[a,B_m(t_m)]$, which, in turn is dominated
by $2^{-m} (g'_+(M) - g'_-(-M))$ for any $t \in [0,K]$.

The second term on the right side of (\ref{eq:trapez}) converges
to 0 as well, since $B_m(t_m)$ converges uniformly to $W(t)$ for
$t \in [0,K]$, and $g'_-$ is bounded on $\mathbb{R}$.

Since the slope of $g$ at any point is bounded above by $g'_+(M)$
and below by $g'_-(-M)$, we see that $g$ is Lipschitz, so
absolutely continuous. Hence
\[
\int_{W(0)}^{W(t)} g'_-(x) \di x = g(W(t)) - g(W(0)) .
\]
Thus we have proved that
\begin{equation} \label{eq:conv1}
\sup_{t \in [0,K]} \left| g(W(t)) - g(W(0)) - T_{x=a}^{B_m(\lfloor
t 2^{2m} \rfloor 2^{-2m})} g'_-(x) 2^{-m} \right| \to 0
\end{equation}
almost surely as $m \to \infty$.

Now we turn to the convergence of the second term on the right
side of (\ref{eq:disc_Ito_Tan_Skor}). Fixing an $\omega \in \Omega
_0$, we want to show that as $m \to \infty$,
\begin{equation} \label{eq:conv2}
\sup_{t \in [0,K]} \left| \sum_{x \in a + 2^{-m} \mathbb{Z}}
\loc^+_m(t,x) (g'_-(x+2^{-m}) - g'_-(x)) - \int_{-\infty}^{\infty}
\frac12 \loc(t,x) \mu (\di x) \right| \to 0
\end{equation}
(The other half containing the $\loc^-_m(t,x)$ terms is similar.)
Now, we have that
\begin{eqnarray} \label{eq:conv_loc_m}
\lefteqn{\sup_{t \in [0,K]} \left| \sum_{x \in a + 2^{-m}
\mathbb{Z}} \loc^+_m(t,x) (g'_-(x+2^{-m}) - g'_-(x)) -
\int_{-\infty}^{\infty} \frac12 \loc(t,x) \mu (\di x) \right|}
\nonumber \\
&\le& \sup_{t \in [0,K]} \sum_{x \in a + 2^{-m} \mathbb{Z}} \left|
\loc^+_m(t,x) - \frac12 \loc(t,x) \right| (g'_-(x+2^{-m})
- g'_-(x))  \\
&+& \frac12 \sup_{t \in [0,K]} \left| \sum_{x \in a + 2^{-m}
\mathbb{Z}} \loc(t,x) (g'_-(x+2^{-m}) - g'_-(x)) -
\int_{-\infty}^{\infty} \loc(t,x) \mu (\di x) \right| \nonumber .
\end{eqnarray}

The first term on the right side of (\ref{eq:conv_loc_m}) goes to
zero, since, by Corollary \ref{co:loctime+}, $\loc^+_m(t,x)$
uniformly converges to $\frac12 \loc(t,x)$ on $[0,K] \times
\mathbb{R}$, and
\[
\sum_{x \in a + 2^{-m} \mathbb{Z}} (g'_-(x+2^{-m}) - g'_-(x)) \le
g'_+(M) - g'_-(-M) < \infty .
\]

The second term on the right side of (\ref{eq:conv_loc_m}) also
goes to zero, since it is the difference of a Riemann--Stieltjes
integral and its approximating sum. Recall that $\loc(t,x)$ is
continuous, non-decreasing in $t$ for any $x$, $\loc(K,x)$ has
compact support, so bounded, as $x$ varies over $\mathbb{R}$, and
the total mass $\mu (\mathbb{R}) < \infty$.

Thus by (\ref{eq:conv1}) and (\ref{eq:conv2}), we have established
the statements of the theorem for measures $\mu$ supported in a
compact interval $[-M,M]$.

Consider now the case of a general $g$. Let $T^M := \inf\lb t
> 0 : |W(t)| \ge M \rb$ for $M=1,2,\dots$, and approximate $g$ by
\[
g^M(x)= \left\{ \begin{array}{ll}
 g(M) + g'_+(M) (x-M) & \mbox{ if } x > M \\
 g(x) & \mbox{ if } |x| \le M \\
 g(-M) + g'_-(-M) (x-M) & \mbox{ if } x < -M .
\end{array} \right.
\]
Define $\mu ^M$ as the measure determined by the change of
$(g^M)'_-$, which is supported in $[-M,M]$ with finite mass.
Clearly, $g(W(t)) = g^M(W(t))$ if $0 \le t \le T^M$ and, since
$\loc(T^M,x) = 0$ for all $x$, $|x| > M$, $\int \loc(t,x) \mu(\di
x) = \int \loc(t,x) \mu^M(\di x)$ if $0 \le t \le T^M$.

The previous argument proves (\ref{eq:non-cadlag_int}) for the
interval $[0,K \wedge T^M]$, and also (\ref{eq:Ito_Tan_form}) if
$0 \le t \le T^M$. Since the stopping times $(T^M)_{M=1}^{\infty}$
increase to $\infty$ a.s., this proves the theorem.

\end{proof}

Comparing It\^o's formula (\ref{eq:Ito_form}) and It\^o--Tanaka
formula (\ref{eq:Ito_Tan_form}) when $g$ is $C^2$ and convex, we
obtain that
\[
\int_{0}^{t} g''(W(s)) \di s = \int_{-\infty }^{\infty } \loc(t,x)
g''(x) \di x .
\]
Since this holds for any continuous and positive function $g''$,
by a monotone class argument we obtain the well-known
\emph{occupation time formula} for any bounded, Borel measurable
function $h$:
\begin{equation} \label{eq:occup}
\int_{0}^{t} h(W(s)) \di s = \int_{-\infty }^{\infty } h(x)
\loc(t,x) \di x .
\end{equation}

As a special case, one gets $\int_{-\infty }^{\infty } \loc(t,x)
\di x  = t$. Also, from (\ref{eq:Ito_Tan_form}) with $g(x) = |x -
a|$, one obtains \emph{Tanaka's formula}
\[
\loc(t,a) = |W(t)-a| - |W(0) - a| - \int_{0}^{t} \sign(W(s)-a) \di
W(s).
\]

%%%%%%%%%%%%%%%%%%%%%%%%%%%%%%%%%%%%%%%%%%%%%%%%%%%%%%%%%%%%%%%%%%%

\section{Integration of predictable processes} \label{sec:Pred}

In this section our aim is to show that when the integrand
$(Y(t))_{t \ge 0}$ is a predictable process satisfying some
integrability condition, our previous approach to approximate the
stochastic integral by sums of values of the integrand at the
Skorohod stopping times $s_m(k)$, multiplied by $\pm 2^{-m}$,
where the signs are obtained by a fair coin tossing sequence,
still works. In other words, the standard non-random partitions of
a time interval $[0,t]$ may be replaced by our random partitions
$0 = s_m(0) < s_m(1) < \cdots < s_m(\lfloor t 2^{2m} \rfloor)$ of
Skorohod embedding of random walks into the Wiener process. This
is not surprising, because, as Lemma \ref{le:equid} shows, such
Skorohod partitions are asymptotically equivalent to partitions
$(k2^{-2m})_{k=0}^{t2^{2m}}$, as $m \to \infty$. This approach
will be described rather briefly, since in this case the details
are essentially the same as in the standard approach to stochastic
integration found in textbooks.

Let $(\mathcal{F}_t)_{t \ge 0}$ denote the natural filtration
defined by our BM $W$. In this paper, stopping times and adapted
processes are understood with respect to this filtration. By
definition, we say that $Y$ is a \emph{simple, adapted process} if
there exist stopping times $0 \le \tau_0 \le \tau_1 \le \cdots$
increasing to $\infty$ a.s. and random variables $\xi _0, \xi _1,
\dots$ such that $\xi _j$ is $\mathcal{F}_{\tau_j}$-measurable,
$\mathbb{E}(\xi _j^2) < \infty$ for all $j$, and
\[
Y(t) = \xi_0 \ind{\{0\}}(t) + \sum_{r=1}^{\infty} \xi_{r-1}
\ind{(\tau_{r-1}, \tau_{r}]}(t) \qquad (t \ge 0).
\]

In the sequel, the only case we will consider is when, with a
given $m \ge 0$, the stopping time sequence is the one given by
the Skorohod embedding (\ref{eq:Skor1}). Let $Y$ be a
left-continuous, adapted process and with a fixed $b > 0$, take
the truncated process
\begin{equation} \label{eq:trunc}
Y^b(t) := \left\{ \begin{array}{rl}
b & \mbox{ if } Y(t) \ge b \\
Y(t) & \mbox{ if } |Y(t)| < b \\
-b & \mbox{ if } Y(t) \le -b .
\end{array} \right.
\end{equation}
Further, with $m \ge 0$ fixed, take
\begin{equation} \label{eq:Ybm}
Y^b_m(t) := Y^b(0) \ind{\{0\}}(t) + \sum_{r=1}^{\infty}
Y^b(s_m(r-1)) \ind{(s_m(r-1), s_m(r)]}(t) \qquad (t \ge 0),
\end{equation}
with the Skorohod stopping times (\ref{eq:Skor1}). Then $Y^b_m$ is
a simple, adapted process.

Then, similarly as in (\ref{eq:stocsum1}), a stochastic sum of
$Y^b_m$ is defined as
\begin{eqnarray} \label{eq:stocsum2}
\left(Y^b_m \cdot W \right)_t &:=& \sum_{r=1}^{\lfloor t 2^{2m}
\rfloor} Y^b(s_m(r-1)) \left(W(s_m(r)) - W(s_m(r-1)) \right)
\nonumber \\
&=& \sum_{r=1}^{\lfloor t 2^{2m} \rfloor} Y^b(s_m(r-1)) X_m(r)
2^{-m},
\end{eqnarray}
where $(X_m(r))_{r=1}^{\infty}$ is a sequence of independent,
$\PP\lb X_m(r)=\pm 1 \rb =\half$ random variables. Observe that
(\ref{eq:stocsum2}) without the $X_m(r) = \pm1$ factors would be
asymptotically an ordinary Riemann sum of the integral
$\int_{0}^{t} Y^b(s) \di s$. (The differences $s_m(r)-s_m(r-1)$
asymptotically equal $2^{-m}$ by Lemma \ref{le:equid}.) The random
$\pm 1$ factors multiplying the terms of the sum transform it into
an approximation of a stochastic integral.

It is clear that the usual properties of stochastic integrals hold
for such stochastic sums. Namely, one has linearity: $(\alpha
Y^b_m + \beta Z^b_m) \cdot W = \alpha Y^b_m \cdot W + \beta Z^b_m
\cdot W$, zero expectation: $\mathbb{E}(Y^b_m \cdot W) = 0$,
isometry: $\mathbb{E}(Y^b_m \cdot W)_K^2 = \int_{0}^{K} \mathbb{E}
(Y^b_m(t))^2 \di t$, etc.

\begin{lem} \label{le:approx}
Let $K > 0$ be fixed and $Y$ be a left-continuous, adapted process
such that
\[
\left\| Y \right\|^2_K := \int_{0}^{K} \mathbb{E} Y^2(t) \di t  <
\infty.
\]
Then there exists a sequence $(m(b))_{b=1}^{\infty}$ such that
$\|Y -Y^b_{m(b)}\|_K \to 0$ as $b \to \infty$, where $Y^b_m$is
defined by (\ref{eq:Ybm}).

\end{lem}

\begin{proof}
First, for any $b$,
\begin{eqnarray*}
\|Y - Y^b\|^2_K &\le& \int_{0}^{K} \mathbb{E} \left(Y^2(t)
\ind{\{|Y(t)| \ge b\}} \right) \di t   \\
&\le& \int_{0}^{K} \mathbb{E} Y^2(t) \di t < \infty.
\end{eqnarray*}
Therefore $\|Y - Y^b\|_K \to 0$ as $b \to \infty$ by dominated
convergence.

Further, for $b$ fixed and $m \to \infty$, $Y^b_m(t) \to Y^b(t)$
a.s. for all $t \in [0,K]$, hence $\|Y^b - Y^b_m\|_K \to 0$ by
bounded convergence. These prove the lemma.

\end{proof}

So let $Y$ be a left-continuous, adapted process such that $\| Y
\|_K < \infty$ and take $J_b(t) := \left(Y^b_{m(b)} \cdot W
\right)_t$, $t \in [0,K]$, where $Y^b_{m(b)}$ is defined by the
previous lemma. Since $Y^b_{m(b)}$ tends to $Y$ in $L^2([0,K]
\times \Omega)$, by isometry, $J_b(t)$ tends to a random variable
$J(t)$ in $L^2(\Omega )$, which is defined as the stochastic
integral $\int_{0}^{t} Y(s) \di W(s)$ for $t \in [0,K]$.

From this point, the discussion of the usual properties of
stochastic integrals and extensions to more general integrands
largely agrees with the standard textbook case, therefore omitted.

%%%%%%%%%%%%%%%%%%%%%%%%%%%%%%%%%%%%%%%%%%%%%%%%%%%%%%%%%%%%%%%%%%%

\section{Extension to continuous local martingales} \label{sec:Mart}

The extension of the methods of the previous sections to
continuous local martingales is rather straightforward. By a
useful theorem of Dambis and Dubins--Schwarz (DDS), a continuous
local martingale $M$ can be transformed into a Brownian motion $W$
by a change of time, when time is measured by the quadratic
variation $\qv{t}$.

Let a right-continuous and complete filtration $(\mathcal{F}_t)_{t
\ge 0}$ be given in $(\Omega , \mathcal{F}, \PP)$ and let the
continuous local martingale $M$ be adapted to it. Define
\[
T_s=\inf \{t \ge 0: \qv{t} \ge s \}  \qquad (s \ge 0).
\]

\begin{thmA} \label{th:DDS1}
If $M$ is a continuous $(\mathcal{F}_t)$-local martingale
vanishing at 0 and such that $\qv{\infty}=\infty$ a.s., then
$W(s)=M(T_s)$ is an $(\mathcal{F}_{T_s})$-Brownian motion and
$M(t)=W(\qv{t})$.

\end{thmA}
Similar statement is true when $\qv{\infty} < \infty$ is possible,
cf. \cite[p. 182]{RevYor1999}. For example, $M(t)=W(\qv{t})$ still
holds for any $t \ge 0$.

Now some results are recalled from \cite{SzeSza2004}.
Skorohod-type stopping times can be defined for $M$, similarly as
for $W$ in (\ref{eq:Skor1}). For $m \ge 0$, let $\tau_m(0)=0$ and
\[
\tau_m(k+1)=\inf \left\{t: t > \tau_m(k), \left| M(t) -
M(\tau_m(k)) \right| = 2^{-m} \right\} \qquad (k \ge 0) .
\]
Then $s_m(k) = \qv{\tau_m(k)}$ and $B_m(k2^{-2m}) = W(s_m(k)) =
M(\tau_m(k))$, where $s_m(k)$ is defined by (\ref{eq:Skor1}) for
the DDS BM $W$ and $B_m$ is a scaled simple, symmetric RW, as in
previous sections.

Let $N_m(t)$ denote the following \emph{discrete quadratic
variation process} of $M$:
\begin{eqnarray*}
N_m(t) &=& 2^{-2m} \# \{r: r>0, \tau_m(r) \le t \}   \\
&=& 2^{-2m} \# \{r: r>0, s_m(r) \le \qv{t} \} \qquad (t\ge 0) .
\end{eqnarray*}
The next theorem is a reformulation of \cite[Theorems 2 and
4]{SzeSza2004}, correcting some unfortunate errors there.

\begin{thmA}
Let $K > 0$ and take a sequence $(c_m)$ increasing to $\infty$
arbitrary slowly. Then for a.e. $\omega  \in \Omega $, if $m \ge
m_0(\omega)$, we have
\[
\sup_{0\le t\le K} \left|\qv{t} - N_m(t)\right| < c_m m^{\half}
2^{-m}
\]
and
\[
\sup_{0\le t\le K} \left| M(t) - B_m(N_m(t)) \right| < c_m
m^{\frac34} 2^{-\frac{m}{2}} .
\]

\end{thmA}

Thus $N_m(t)$ almost surely uniformly converges on compacts to
$\qv{t}$, and $B_m(N_m(t))$ similarly converges to $M(t)$.

Combining the DDS time change and the pathwise approximations of
local time described in Section \ref{sec:Pre}, it is not too hard
to generalize the local time results of \cite{SzaSze2005} to
continuous local martingales $M$. Let the $m$th approximate local
time of $M$ at a point $x \in 2^{-m} \mathbb{Z}$ up to time $t \in
2^{-2m} \mathbb{Z}_+$ be defined by
\[
\loc^{M,\pm}_m(t,x) = 2^{-m} \#\{s : \: B_m(s) = x, \:
B_m(s+2^{-2m}) = x \pm 2^{-m} \} ,
\]
where $s \in 2^{-2m} \mathbb{Z} \cap \left[0,\qv{t}\right)$ and
$B_m(j 2^{-2m}) = M(\tau _m(j))$. Then $\loc_m^{M,\pm}(t,x)$ can
be extended to any $t \in \mathbb{R}_+$ and $x \in \mathbb{R}$ by
linear interpolation, as a continuous two-parameter process.

The major difference compared to Corollary \ref{co:loctime+} is
that that the time interval $[0,K]$ has to be replaced by an
interval of the form $[0, T^K_b]$ here, where $T^K_b := T_b \wedge
K$, over which $\qv{t}$ is bounded by a constant $b$.

\begin{cor} \label{co:mart_loctime+}
On any strip $[0, T^K_b] \times \mathbb{R}$, the approximate local
times $(\loc^{M,+}_m(t,x))$ a.s. uniformly converge as $m \to
\infty$ to the one half of the local time $\loc^M(t,x) =
\loc^W(\qv{t},x)$, where $\loc^W$ denotes the local time of the
DDS BM of $M$. For any $C>1$, and for any $K>0$, $b \ge e$, and $m
\ge 1$ such that $b 2^{2m} \ge N(C)$, we have
\begin{eqnarray*}
\PP \lb \sup_{(t,x) \in [0,T^K_b] \times \mathbb{R}} \left|\half
\: \loc^M(t,x) - \loc^{M,+}_{m+1}(t,x)\right| \ge 50 \: C
b^{\quart} (\log b)^{\frac34} m^{\frac34} 2^{-\frac{m}{2}} \rb \\
\le \frac{32}{1-4^{1-C}} (b 2^{2m})^{1-C} .
\end{eqnarray*}
Similar statements hold for $(\loc^-_m(t,x))$ as well.
\end{cor}

From this, by the Borel--Cantelli lemma, we get that for a.e.
$\omega \in \Omega $, if $m \ge m_0(\omega)$, we have
\begin{eqnarray*}
\lefteqn{\sup_{(t,x) \in [0,T^K_b] \times \mathbb{R}} \left|\half
\: \loc^M(t,x) -
\loc^{M,+}_{m+1}(t,x)\right|} \\
&=& \sup_{(t,x) \in [0,K] \times \mathbb{R}} \left|\half \:
\loc^M(t \wedge T^K_b,x) - \loc^{M,+}_{m+1}(t \wedge
T^K_b,x)\right| \\
&\le& 50 \: C b^{\quart} (\log b)^{\frac34} m^{\frac34}
2^{-\frac{m}{2}} .
\end{eqnarray*}

If we replace $b$ here by a sequence $(c_m)$ increasing to
$\infty$, then, by the continuity of $\qv{t}$, $T_{c_m} \nearrow
\infty$, and so $t \wedge T^K_{c_m} \nearrow t$ for $t \in [0,K]$.
Hence we have
\begin{cor} \label{co:mart_loctime}
Let $K > 0$ and take a sequence $(c_m)$ increasing to $\infty$
arbitrary slowly. Then for a.e. $\omega  \in \Omega $, if $m \ge
m_0(\omega)$, we have
\[
\sup_{(t,x) \in [0,K] \times \mathbb{R}} \left|\half \:
\loc^M(t,x) - \loc^{M,\pm}_{m}(t,x)\right| < c_m m^{\frac34}
2^{-\frac{m}{2}} .
\]

\end{cor}

Next define It\^o integrals w.r.t. $M$ when the integrand is a
`nice' function of $M$. An It\^o sum is defined by
\begin{eqnarray} \label{eq:stocsum3}
\left(f(M) \cdot M \right)^m_t &:=& \! \! \! \sum_{\tau_m(r) \le
t} f(M(\tau_m(r-1))) \: 2^{-m}X_m(r)
\nonumber \\
&=& \! \! \! \! \! \! \sum_{s_m(r) \le \qv{t}} f(W(s_m(r-1))) \:
2^{-m}X_m(r) ,
\end{eqnarray}
where $2^{-m}X_m(r) = M(\tau_m(r)) - M(\tau_m(r-1)) = W(s_m(r)) -
W(s_m(r-1))$, $(X_m(r))_{r=1}^{\infty}$ being an independent, $\pm
1$ symmetric coin tossing sequence.

\begin{thm} \label{th:mart_non-cadlag_int}
Let $g$ be a $C^2$-function, alternatively, be the difference of
two convex functions. Then for arbitrary $K>0$,
\[
\sup_{t\in[0,K]}\left|\left(g'_- (M) \cdot M \right)_t^m -
\int_0^t g'_- (M(s)) \di M(s) \right| \rightarrow 0
\]
almost surely as $m \to \infty$, and in the first case one gets
the It\^o formula
\[
g(M(t)) - g(M(0)) = \int_{0}^{t} g'(M(s)) \di M(s) + \frac12
\int_{0}^{t} g''(M(s)) \di \qv{s} ;
\]
in the second case one obtains the It\^o--Tanaka formula
\[
g(M(t)) - g(M(0)) = \int_{0}^{t} g'_- (M(s)) \di M(s) + \frac12
\int_{-\infty}^{\infty} \loc^M (t,x) \mu(\di x)  .
\]

\end{thm}

The proof is similar to the Brownian case, therefore omitted.

(\ref{eq:stocsum3}) indicates that if the function $f$ is $C^1$ or
of finite variation on bounded intervals, then
\[
\int_0^t f(M(s)) \di M(s) = \int_{0}^{\qv{t}} f(W(s)) \di W(s).
\]

We mention that the DDS BM $W$ and the quadratic variation $\qv{}$
are independent processes if and only if $M$ is symmetrically
evolving (is an Ocone martingale) in the sense that given the past
of $M$, distributions of future increments are invariant under
reflection, cf. \cite[Theorem 5]{SzeSza2004}.

Finally, we can extend the method of integration of predictable
integrands discussed in Section \ref{sec:Pred} to integrals w.r.t.
a continuous local martingale $M$. Let $K > 0$ be fixed and $Y$ be
a left-continuous, adapted process on $[0, K]$ such that
\[
\left\| Y \right\|^2_{M,K} := \mathbb{E} \int_{0}^{K} Y^2(t) \di
\qv{t} < \infty.
\]
For $b > 0$ let $Y^b$ denote the truncated process
(\ref{eq:trunc}). If $m \ge 0$ fixed, take
\begin{equation} \label{eq:mart_Ybm}
Y^b_m(t) := Y^b(0) \ind{\{0\}}(t) + \sum_{r=1}^{\infty}
Y^b(\tau_m(r-1)) \ind{(\tau_m(r-1), \tau_m(r)]}(t) ,
\end{equation}
for $t \ge 0$. Then $Y^b_m$ is a simple, adapted process.

An It\^o sum of $Y^b_m$ is defined as
\begin{eqnarray*}
\left(Y^b_m \cdot M \right)_t &:=& \sum_{\tau_m(r) \le t}
Y^b(\tau_m(r-1)) \left(M(\tau_m(r)) - M(\tau_m(r-1)) \right)
 \\
&=& \sum_{\tau_m(r) \le t} Y^b(\tau_m(r-1)) \: X_m(r) 2^{-m},
\end{eqnarray*}
where $(X_m(r))_{r=1}^{\infty}$ is a sequence of independent,
$\PP\lb X_m(r)=\pm 1 \rb =\half$ random variables.

Here too, cf. Lemma  \ref{le:approx}, there exists a sequence
$(m(b))_{b=1}^{\infty}$ such that $\|Y -Y^b_{m(b)}\|_{M,K} \to 0$ as
$b \to \infty$, where $Y^b_m$ is defined by (\ref{eq:mart_Ybm}).
So take $J_b(t) := \left(Y^b_{m(b)} \cdot M \right)_t$, $t \in
[0,K]$. Since $Y^b_{m(b)}$ tends to $Y$ in $L^2_M([0,K] \times
\Omega)$, by isometry, $J_b(t)$ tends to a random variable $J(t)$
in $L^2(\Omega )$, which is defined as the stochastic integral
$\int_{0}^{t} Y(s) \di M(s)$ for $t \in [0,K]$.

%%%%%%%%%%%%%%%%%%%%%%%%%%%%%%%%%%%%%%%%%%%%%%%%%%%%%%%%%%%%%%%%%%%

\end{document}